\def\to{\rightarrow}

\def\be{\begin{equation}}

\def\bs{\bigskip}\def\ms{\medskip}
\def\no{\noindent}

\def\const{\text{\rm const}}

\def\ker{\text{\rm ker}}

\def\dist{\text{\rm dist}}
\def\supp{\text{\rm supp}\,}

\def\GG{{{\textsl{G}}}}

\def\TT{{\textsl{T}}}

\def\e{{\epsilon}}
\def\lan{{\lambda_{n}}}
\def\Th{{\Theta}}
\def\Tb{{\bar{\Theta}}}
\def\L{\Lambda}

\def\BM{{\mathcal BM}}

\def\C{{\mathbb{C}}}

\def\Z{{\mathbb{Z}}}
\def\N{{\mathbb{N}}}
\def\R{{\mathbb R}}

\documentstyle[10pt,amssymb]{amsart}

\title{P\'olya sequences, Toeplitz kernels and gap theorems}
\author{Mishko Mitkovski}

\address{Texas A\&M University
\\ Department of Mathematics\\
College Station, TX 77843, USA } \email{mmitkov@@math.tamu.edu}
\author{Alexei Poltoratski}
\address{Texas A\&M University
\\ Department of Mathematics\\
College Station, TX 77843, USA }
\email{alexeip@@math.tamu.edu}
\thanks{The second author is supported by
 N.S.F. Grant No. 0800300}

\theoremstyle{plain}

\newtheorem*{lem1}{Lemma 1}
\newtheorem*{lem2}{Lemma 2}

\newtheorem*{thm1}{Theorem I}
\newtheorem*{thm2}{Theorem II}
\newtheorem*{thm3}{Theorem III}

\newtheorem*{thmA}{Theorem A}
\newtheorem*{thmC}{Theorem C}
\newtheorem*{thmB}{Theorem B}

\newtheorem*{cor}{Corollary}

\theoremstyle{definition}
\newtheorem*{remark}{Remark}

\newenvironment{proof}
{{\noindent\it Proof:}}{\hfill$\Box$}

\newenvironment{proofA}
{{\noindent{\bf \textit{Proof of Theorem A:}}}}{\hfill$\Box$}

\newenvironment{proofB}
{{\noindent\bf \textit{Proof of Theorem B:}}}{\hfill$\Box$}

\newenvironment{proofC}
{{\noindent\bf \textit{Proof of Theorem C:}}}{\hfill$\Box$}

\numberwithin{equation}{section}

\begin{document}
\maketitle
\begin{abstract}
A separated sequence $\Lambda$ on the real line is called a P\'olya sequence if
any entire function of zero exponential type bounded on $\Lambda$ is constant.
In this paper we solve the problem by P\'olya and Levinson that asks for a description
of P\'olya sets. We also show that the P\'olya-Levinson  problem is equivalent to a version
of the so-called Beurling gap problem on  Fourier transforms
of measures. The solution is obtained via a recently developed approach based on the use of Toeplitz kernels
and de Branges spaces of entire functions.
\end{abstract}

\bs\section{Introduction and background}

\bs\subsection{Introduction}
An entire function $F$ is said to have exponential type zero if
$$\limsup_{|z|\to \infty}\frac{\log{|F(z)|}}{|z|}=0.$$
We call a separated real sequence
$\Lambda=\{ \lan \}_{n=-\infty}^{\infty}$ (a sequence satisfying \linebreak
$| \lan - \lambda_{m} | \geq \delta
>0, \hspace{.2cm} (n \neq m)$) a {\it P\'olya sequence} if any entire function
of exponential type zero that is bounded on $\Lambda$ is constant.
In this paper we consider the problem of description of P\'olya sequences.

Historically, first results on P\'olya sequences were obtained in
the work of Valiron~\cite{Val}, where it was proved that the set
of integers $\Z$ is a P\'olya sequence. Later this result was
popularized by P\'olya, who posted it as a problem in~\cite{Pol}.
Subsequently many different proofs and generalizations were given
(see for example section 21.2 of~\cite{Lv1} or chapter 10 of
~\cite{Boa} and references therein).

In his 1940' book \cite{Lev} Levinson showed that if $|\lan-n|\leq
p(n)$, where $p(t)$ satisfies $\int \frac{p(t)}{1+t^2}\log
|\frac{t}{p(t)}|dt < \infty$ and some smoothness conditions, then
$\Lambda=\{ \lan \}$ is a P\'olya sequence. In the same time for
each such $p(t)$ satisfying $\int{p(t)dt/(1+t^2)}=\infty$ he was
able to construct a sequence $\Lambda=\{ \lan \}$ that is not
P\'olya sequence. As it often happens in problems from this area,
the construction took a considerable effort (see \cite{Lev}, pp.
153-185). Closing the gap between Levinson's sufficient condition
and the counterexample remained an open problem for almost 25
years until de Branges \cite{dBr1} essentially solved it by
showing that $\Lambda$ is a P\'olya sequence if
$\int{p(t)dt/(1+t^2)}<\infty$ (but assuming extra regularity
conditions on the sequence).

The results of \cite{Lev} and \cite{dBr1} remain strongest to
date. However, none of them gives a complete answer, since there
are P\'olya sequences for which $\int{p(t)dt/(1+t^2)}=\infty.$ For
example, as will be clear from our results below,  the sequence
$$\lan:=n+n/\log{(|n|+2)}, \hspace{.1cm} n\in \Z$$ is a P\'olya
sequence.

In the opposite direction, \cite{dBr1} contains the following
necessary condition. A sequence of disjoint intervals $I_n$ on the
real line is called long (in the sense of Beurling and Malliavin)
if $$\sum_n\frac {|I_n|^2}{1+\dist^2(I_n,0)}=\infty,$$ and it is
called short otherwise. Here $|I_n|$ denotes the length of the interval
$I_n$. De~Branges \cite{dBr1} proved that if the
complement of a closed set $X \subset \R$ is long then there exists a non
constant zero type entire function that is bounded on $X$. In
particular, if the complement of a sequence $\Lambda$ is long then
$\Lambda$ is not a P\'olya sequence. The sequence $\lan:=n^2$ shows
that this condition is not sufficient. Indeed, $\lan=n^2$ is the
zero set of the zero type function $F(z):=\cos{\sqrt{2\pi
z}}\cos{\sqrt{-2\pi z}}$ and thus is not a P\'olya sequence. On the
other hand, the real complement of this sequence is short.

In this paper we give the following answer to the P\'olya-Levinson
question, see the corollary of Theorem A below. We show that a separated sequence
of real numbers $\Lambda$ is \emph{not} a P\'olya sequence if and only if there exists
a long sequence of intervals $\{I_n \}$ such that
$$\frac{ \#(\Lambda \cap I_n)}{|I_n|}\to 0.$$

Our approach is similar to the one developed by the second author
and N.~Makarov in~\cite{MP} and~\cite{MP2}, where it was used to
obtain extensions and applications of  the Beurling-Malliavin
theory. One of our main tools
is the connection between the P\'olya-Levinson problem, the gap problem
and the problem of injectivity of Toeplitz operators, see Theorems
A and C. To promote the Toeplitz approach and to make the paper more self-contained
we often include full proofs rather than referring to existing results.

The Beurling gap problem that we consider here may be formulated
as follows. Under what conditions on a separated real sequence
$\Lambda=\{ \lan \}_{n=-\infty}^{\infty}$ does there exist a
nonzero finite measure $\mu$ supported on $\Lambda$ such that the
Fourier transform of $\mu$ vanishes on an interval of positive
length? Of course, one can ask the same question for an arbitrary
closed set $X \subset \R$. It is interesting to observe that the
sufficient condition that de Branges gave for the P\'olya-Levinson
problem was the same as Beurling's sufficient condition for the
gap problem, see \cite{Beu} or \cite{Koo1}. Namely, Beurling
proved that if the complement of a closed set $X \subset \R$ is
long then there exists no nonzero finite measure $\mu$ supported
on $X$ such that the Fourier transform of $\mu$ vanishes on an
interval of positive length. As the reader will see below, this
was not a coincidence since the problems for sequences turn out to
be equivalent.

We give a solution to the gap problem for separated sequences and
 improve Beurling's gap theorem for
general closed sets, see Theorem B and its corollaries.

\ms{\bf Acknowledgment.}
We are grateful to N.~Makarov who brought the P\'olya-Levinson problem to our attention
and to M.~Sodin for useful discussions and references.

\subsection{Background} We use the standard notation $N^+(\C_+)$
to denote the Smirnov-Nevanlinna class in the upper half-plane
$\C_+=\{z| \Im{z}>0\}$ consisting of analytic functions $f(z)$
that can be represented as a ratio $g(z)/h(z)$ of two bounded
analytic functions with $h(z)$ being outer. Each function in
$N^+(\C_+)$ has non-tangential boundary values almost everywhere
on $\R$ that completely determine the function.  A mean type of a
function $f(z)$ in $N^+(\C_+)$ is defined as $\tau:=\limsup_{y \to
\infty}\log{|F(iy)|}/y$. It is easy to see that every function in
$N^+(\C_+)$ has a non-positive mean type which is exactly the
exponent of $S(z):=e^{iz}$ in the inner-outer factorization of
$f(z)$ taken with a negative sign. Here and throughout the paper
$S(z)$ denotes the singular inner function $e^{iz}$.

A Hardy space $H^2(\C_+)$ consists exactly of those functions in
$N^+(\C_+)$ which are square-integrable on $\R$ (for more on
Smirnov-Nevanlinna and Hardy spaces see, e.g.,~\cite{Ga}).

A classical theorem of Krein gives a connection between the
Smirnov-Nevanlinna class $N^+(\C_+)$ and the Cartwright class
$C_a$ consisting of all entire functions $F(z)$ of exponential
type $\leq a$ that satisfy $\log|F(t)| \in L^1(dt/(1+t^2))$. An
entire function $F(z)$ belongs to the Cartwright class $C_a$ if
and only if
$$ \frac{F(z)}{S^a(z)} \in N^+(\C_+), \hspace{1cm}
\frac{F^{\#}(z)}{S^a(z)} \in N^+(\C_+),$$ where
$F^\#(z)=\overline{F(\bar{z})}$.

As an immediate consequence one obtains a connection between the
Hardy space $H^2(\C_+)$ and the Paley-Wiener space $PW_a$. Namely,
an entire function $F(z)$ belongs to the Paley-Wiener class $PW_a$
if and only if
$$ \frac{F(z)}{S^a(z)} \in H^2(\C_+), \hspace{1cm}
\frac{F^{\#}(z)}{S^a(z)} \in H^2(\C_+).$$

The definition of the de~Branges spaces of entire functions may be
viewed as a generalization of the above definition of the
Payley-Wiener spaces with $S^a(z)$ replaced by a more general
entire function. Consider  an entire function  $E(z)$ satisfying
the inequality
$$ |E(z)|>|E(\bar{z})|, \hspace{1cm} z \in \C_{+}.$$ Such
functions are usually called  de~Branges functions. The de~Branges
space $B_E$ associated with $E(z)$ is defined to be the space of
entire functions $F(z)$ satisfying
$$ \frac{F(z)}{E(z)} \in H^2(\C_+), \hspace{1cm}
\frac{F^{\#}(z)}{E(z)} \in H^2(\C_+).$$ It is a Hilbert space
equipped with the norm $\|F\|_E:=\|F/E\|_{L^2(\R)}.$ By a
theorem of Krein, if $E(z)$ is of exponential type then
all the functions in the de~Branges space $B_E$ will be of
exponential type not greater then the type of $E(z).$ A de~Branges
space is called short (or regular) if together with
every function $F(z)$ it contains $(F(z)-F(a))/(z-a)$ for any
$a\in\C$.

We will utilize the following  well-known result from the theory
of de Branges spaces of entire functions.

\begin{thm1}\cite{dBr2}
Let $\mu$ be a positive measure on $\R$ satisfying
$\int{d\mu(t)/(1+t^2)}<\infty.$ Then there exists a short
de~Branges space $B_E$ contained isometrically in $L^2(\mu)$, with
de~Branges function $E(z)$ being of Cartwright class and having no
real zeros. Moreover, if there exists such a space $B_E$ with  $E(z)$
 of positive exponential type, then there also exists such a  space $B_E$ that is contained properly in $L^2(\mu)$.
\end{thm1}

\begin{remark} The existence part follows from Theorem XII of
\cite{dBr2}. The second part follows from Theorems IV and X in the
same paper. Finally, the shortness can be derived from the proof of Theorem XII \cite{dBr2} or
from problem 71
in~\cite{dBr} by taking $S(z)=1$.
\end{remark}

General treatment of de~Branges' theory is given in~\cite{dBr}.

\ms

Every de~Branges function gives rise to a meromorphic inner
function $\Th(z)=E^{\#}(z)/E(z).$ We say that an inner function
$\Th(z)$ in $\C_+$ is a meromorphic inner function if it allows a
meromorphic extension to the whole complex plane. The meromorphic
extension to the lower half-plane $\C_{-}$ is given by:
$$\Th(z)=\frac{1}{\Th^{\#}(z)}. $$ Conversely, every
meromorphic inner function $\Th(z)$ can be represented in the form
$\Th(z)=E^{\#}(z)/E(z),$ for some de~Branges function $E(z)$ (see,
for instance, Lemma 2.1 in~\cite{MH}). Such a function is unique
up to a factor of an entire function that is real on $\R$ and has only real zeros. 

Each inner function $\Th(z)$ determines a model subspace
$$K_\Theta=H^2\ominus \Theta H^2$$ of the Hardy space
$H^2(\C_+)$. These subspaces play an important role in complex and
harmonic analysis, as well as in operator theory,
see~\cite{Ni2,Ni}. There is an important relationship between the
model subspaces $K_\Th$ and the de~Branges spaces $B_E$ of entire
functions. If $E(z)$ is a de~Branges function and
$\Th(z)=E^{\#}(z)/E(z)$ is the corresponding meromorphic inner
function, then the multiplication operator
$f\mapsto Ef$ is an isometric isomorphism $K_\Th \rightarrow
B_E$.

Each inner function $\Th(z)$ determines a positive harmonic
function $\Re \frac{1+\Th(z)}{1-\Th(z)}$ and by a Herglotz
representation a positive measure $\sigma$ such that
\begin{equation} \label{for1} \Re
\frac{1+\Th(z)}{1-\Th(z)}=py+\frac{1}{\pi}\int{\frac{yd\sigma
(t)}{(x-t)^2+y^2}}, \hspace{1cm} z=x+iy,\end{equation} for some $p
\geq 0$. The number $p$ can be viewed as a point mass at infinity.
The measure $\sigma$ is singular, supported on $\{ \Th=1 \}
\subset \R,$ and satisfies $\int{d\sigma(t)/(1+t^2)}<\infty.$ It
is usually called a Clark measure for $\Th(z)$. Conversely, for
every positive singular measure $\sigma$ with
$\int{d\sigma(t)/(1+t^2)}<\infty$ and a number $p \geq
0$, there exists an inner function $\Th(z)$ determined by the
formula~\eqref{for1}. Below, when we say that an inner function
$\Th(z)$ corresponds to $\sigma$ we always assume $p=0$.

Every function $f \in K_\Th$ can be represented by the formula
\begin{equation} \label{for2} f(z)=\frac{p}{2\pi
i}(1-\Th(z))\int{f(t)\overline{(1-\Th(t))}dt}+\frac{1-\Th(z)}{2\pi
i}\int{\frac{f(t)}{t-z} d\sigma (t)}.
\end{equation}
If $1-\Th(t) \notin L^2(\R)$ then $p=0$ and hence we have a nicer
looking formula
$$f(z)=\frac{1-\Th(z)}{2\pi
i}\int{\frac{f(t)}{t-z} d\sigma (t)}.$$ This gives an isometry of
$L^2(\sigma)$ onto $K_\Th$. In the case of meromorphic $\Th(z)$,
every function $f \in K_\Th$ also has a meromorphic extension in
$\C$, and it is given by the formula~\eqref{for2}. The
corresponding Clark measure is discrete with atoms at the points
of $\{\Theta=1\}$ given by $\sigma(\{x\})=\frac{2\pi
}{|\Th'(x)|}$.

Each meromorphic inner function $\Th(z)$ can be written as
$\Th(t)=e^{i\phi(t)}$ on $\R$, where $\phi(t)$ is a real analytic
and strictly increasing function. The function $\phi(t)=\arg{\Th(t)}$ is
the continuous argument of $\Th(z)$. The phase function of $E(z)$ is defined as
$-\frac{1}{2}\arg{\Th(t)},$ where $\Th(z)$ is the corresponding
meromorphic inner function.

A subset of $\R$ is called discrete if it has no finite density points.
For every discrete set  $\Lambda\subset\R$, there exists a (far from unique) meromorphic inner
function $\Th(z)$ such that $\{\Th=1\}=\Lambda.$ In the case of a
separated sequence $\Lambda$, there is a meromorphic inner
function $\Th(z)$ with $\{\Th=1\}=\Lambda$ whose continuous
argument $\arg{\Th(t)}$ has a bounded derivative (see for instance Lemma~16
in~\cite{dBr2}).

Recall that the Toeplitz operator $T_U$ with a symbol $U\in
L^\infty(\R)$ is the map
$$T_U:H^2\to H^2,\qquad F\mapsto P_+(UF),$$
where $P_+$ is the orthogonal projection  in $L^2(\R)$ onto the
Hardy space $H^2=H^2(\C_+)$.

We will use the following notation for kernels of Toeplitz
operators (or {\it Toeplitz kernels} in $H^2$):
$$N[U]=\ker{T_U}.$$
For example, $N[\bar\Theta]=K_\Theta$ if $\Theta$ is an inner
function. Along with $H^2$-kernels, one defines Toeplitz kernels
in the Smirnov class $N^+(\C_+)$, $$ N^+[U]=\{f \in N^+ \cap
L^1_{loc}(\R): \bar{U}\bar{f} \in N^+\}.$$

\subsection{Beurling-Malliavin densities} \bs\no Before we formulate our results
let us discuss the following notion of density of a discrete sequence and related theorems.

Following \cite{BM2} we say that a discrete sequence $\L\subset
\R$ is $a$-\textit{regular} if for every $\epsilon>0$ any sequence of
disjoint intervals $\{I_n \}$ that satisfies
$$\left|\frac{\#(\L\cap I_n)}{|I_n|}-a\right|\geq
\epsilon$$
for all $n$, is short.

A slightly different $a$-regularity can be defined in the following way,
that is more convenient in some settings.
For a discrete sequence $\L\subset \R$ we denote by $n_\L (x)$ its
continuous counting function, i.e. the function that is
continuous on $\R$, grows linearly by $1$ between each pair of neighboring
points of $\L$ and is equal to $0$ at $0$.
We say that $\L$ is \textit{strongly} $a$-\textit{regular}
if

$$\int\frac{|n_{\L}-ax|}{1+x^2}<\infty.$$

\no Conditions like this can be found in many related results, see for instance
\cite{dBr} or \cite{Koo1}. Even though $a$-regularity is not equivalent to
strong $a$-regularity, in the following definitions of densities changing
"$a$-regular" to "strongly $a$-regular"  will lead to equivalent definitions.

\bs

The interior BM (Beurling-Malliavin)
density is defined as

\begin{equation}\label{id}D_*(\L):=\sup \{a\ |\ \exists \ \textrm{$a$-regular subsequence}\ \ \L'\subset \L \}.\end{equation}

Similarly, the exterior BM density is defined as

\begin{equation}\label{ed}D^*(\L):=\inf \{a\ |\ \exists \ \textrm{$a$-regular supsequence}\ \ \L'\supset \L \}.\end{equation}
\no If no such sequence exists $D^*(\L):=\infty$, see \cite{BM2}. Exterior density was
used in the Beurling-Malliavin solution of the completeness problem for families
of exponential functions in $L^2$ on an interval, see \cite{BM2}, \cite{HJ} or \cite{Koo1}.

The following simple observation will be useful for us in the next
section: $D_*(\L)=0$ if and only if there exists a long sequence
of intervals $\{I_n\}$ such that

$$\#(\L\cap I_n)=o(|I_n|)\ \ as \ \ |n|\to\infty.$$

\bs

A description of $D^*(\L)$ in terms of Toeplitz kernels is given
by the following formula, see ~\cite{MP}:
$$D^*(\L)=\frac{1}{2\pi}\sup\{a: N[\bar{S}^a\Th]=0\}, $$ where $\Th(z)$
denotes some/any meromorphic inner function with
$\{\Th=1\}=\Lambda.$ \no Below (see Theorems B and C) we give a
similar description of the interior BM density for separated
sequences $\L$. Namely,
$$D_*(\L)=\frac{1}{2\pi}\sup\{a: N[\Tb S^a]=0\}, $$ where $\Th(z)$
denotes some/any meromorphic inner function with
$\{\Th=1\}=\Lambda.$

An equivalent way to define the interior BM density is as follows.
 Let $\gamma: \R\to\R$
be a continuous function such that $\gamma(\mp\infty)=\pm\infty$.
i.e.
$$\lim_{x\to-\infty}\gamma(x)=+\infty,\qquad \lim_{x\to+\infty}\gamma(x)=-\infty.$$
The family $\BM(\gamma)$ is defined as the collection  of the connected
components of the open set
$$\left\{x \in \R:~\gamma(x)\ne\max_{t \in [x,+\infty)}\gamma(t) \right\}.$$
 We say that
$\gamma$ is almost decreasing if
$\gamma(\mp\infty)=\pm\infty$
and
the family  of the intervals $\BM(\gamma)$ is  short.

Now we can state an equivalent definition for interior
BM~density:

\begin{equation}\label{id2}D_*(\L):=\sup \{a\ |\ ax-n_\L (x) \textrm{ is almost decreasing}\}.\end{equation}
Equivalence of this definition and  \eqref{id} can be easily verified.



We will use the following formulations of the Beurling-Malliavin
theorems~\cite{BM1,BM2}:

\begin{thm2}[\cite{MP}] Suppose that $\Th(z)$ is a meromorphic inner
function with the derivative of $\arg\Th(t)$ bounded on $\R$. Then
for any meromorphic inner function $J(z),$ we have $$
N^+[\bar{\Th}J] \neq 0 \hspace{.5cm} \Rightarrow \hspace{.5cm}
\forall \epsilon
>0, \hspace{.5cm} N[\bar{S}^\epsilon \bar{\Th}J] \neq 0.$$
\end{thm2}

\begin{thm3}[\cite{MP}] Suppose $\gamma'(t) > -\const.$
\begin{itemize}
\item[\textit{(i)}] If $\gamma$ is not almost decreasing, then for
every $\epsilon>0$, $N^+[S^\epsilon e^{i\gamma}]=0$.

\item[\textit{(ii)}] If $\gamma$ is almost decreasing, then for
every
 $\epsilon>0$, $N^+[\bar{S}^\epsilon e^{i\gamma}] \neq 0$.
\end{itemize}
\end{thm3}

\begin{remark} As noted in~\cite{MP}, the part $(i)$ of Theorem~III holds
without the assumption $\gamma'(t) > -\const.$
\end{remark}

\bs\section{Results and proofs}

\bs\subsection{Main results} \bs\no As was mentioned in the
introduction, a sequence of real numbers is called separated if $
| \lan - \lambda_{m} | \geq \delta>0, \hspace{.2cm} (n \neq m) $.
It is natural to introduce a separation condition in the
P\'olya-Levinson problem because of the following obvious reasons.
If one takes a zero set of a zero-type entire function  and adds a
large number of points close enough to each zero, the entire
function will still be bounded on the new sequence. At the same
time, this way one can obtain non-P\'olya sequences of arbitrarily
large density, in any reasonable definition of density. Hence, if
one hopes to obtain a description of P\'olya sequences based on
densities or similar terms,  it is necessary  to include a
separation condition, as it  was done in the classical results
cited above.

Recall that a separated sequence
 $\{
\lan \}_{n=-\infty}^\infty$ is called a \textit{P\'olya sequence} if
every zero-type entire function  bounded on $ \{ \lan \} $ is
constant.

\begin{thmA} \label{thmA} Let $\Lambda=\{ \lan
\}_{n=-\infty}^\infty \subset \R$ be a separated sequence of real
numbers. The following are equivalent:
\begin{itemize}
\item[\textit{(i)}] $\Lambda=\{\lan\}$ is a P\'olya sequence.

\item[\textit{(ii)}] There exists  a non-zero measure $\mu$ of
finite total variation, supported on $ \Lambda $, such that the
Fourier transform of $\mu$ vanishes on an interval of positive length.

\item[\textit{(iii)}] The interior Beurling-Malliavin density of
$\L$, $D_*(\L)$, is positive.

\item[\textit{(iv)}] There exists a meromorphic inner function
$\Th(z)$ with $\{ \Th =1 \} = \{ \lan \} $ such that $ N[\Tb
S^{2c}] \neq 0$, for some $c>0$.
\end{itemize}
\end{thmA}

As an immediate consequence we obtain that the sequence of
integers $\Z$ is a P\'olya sequence, as known from Valiron's
original statement. This follows from $(iii)$ and also from $(iv)$
by taking $\Th(z)=S^{2\pi }(z)$. Another consequence is that a
separated real sequence with density zero cannot be a P\'olya
sequence. However, there are sequences with positive density (and
hence positive exterior Beurling-Malliavin density) which are not
P\'olya. As was mentioned in the introduction, the first example
of such a sequence was given by Levinson~\cite{Lev}. New examples
in both directions can now be constructed using the following
description.

\begin{cor} Let $\Lambda=\{ \lan \}_{n=-\infty}^\infty$ be a separated
sequence of real numbers. Then $\Lambda$ is a P\'olya sequence if
and only if for every long sequence of intervals $\{I_n \}$ the
sequence $\frac{ \#(\Lambda \cap I_n)}{|I_n|}$ is not a null
sequence, i.e., $\frac{ \#(\Lambda \cap I_n)}{|I_n|} \nrightarrow
0$.
\end{cor}

In regard to the gap problem we will prove the following result.
Let $M$ denote the set of all complex measures of finite total
variation on $\R$. For  $\mu\in M$ its Fourier transform $\hat \mu(x)$ is defined
as

$$\hat\mu(x)=\int{e^{ixt}d\mu (t)}.$$
If $X$ is a closed subset of the real line denote by $\GG(X)$ the
gap characteristic of $X$:
$$\GG(X):=\sup\{a\ |\ \exists\ \mu\in M,\ \mu\not\equiv 0,\ \supp\mu\subset X,\ \text{ such that}\ \hat{\mu}=0\ \text{ on }\ [0,a] \}.$$\newpage

\begin{thmB} The following are true:
\begin{itemize}
\item[\textit{(i)}] For any separated sequence $\L\subset\R$,
$\GG(\L)\geq 2\pi D_*(\L)$.

\item[\textit{(ii)}] For any closed set $X\subset\R$, $\GG(X)\leq
2\pi D_*(X)$.
\end{itemize}
\end{thmB}

\ms

\begin{cor}
For separated sequences $\L\subset \R$, $\GG(\L)=2\pi D_*(\L)$.
\end{cor}

\ms

The formula for $\GG(X)$ for a general closed set $X$ is more
involved, see \cite{P}. Another immediate consequence of Theorem B
is the following extension of Beurling's gap theorem:

\ms

\begin{cor} Let $X$ be a closed subset of the real line. If there
exists a long sequence of intervals $\{I_n \}$ such that $$\frac{
\#(X \cap I_n)}{|I_n|} \rightarrow 0$$ then any measure $\mu$ of finite total variation supported on $X$,
whose Fourier transform  vanishes on
an interval of positive length, is trivial.
\end{cor}

\ms

Finally, our next result connects the gap problem to the
problem on injectivity of Toeplitz operators. It provides one
of the main tools for our proofs.

If $X\subset \R$ is closed, define

$$\TT(X)=$$
$$
\sup\{a\ |\ \exists\ \textrm{ meromorphic inner }\Th(z)
\textrm{ with }\{\Th=1\}\subset X\textrm{ and }N[\Tb  S^a]\neq 0\}.$$



\ms

\begin{thmC} For any closed $X\subset \R$,

$$\TT(X)=\GG(X).$$

\end{thmC}

\ms

Theorems A, B and C will be proved in the last section.

\bs\subsection{Technical lemmas}\no For the main proofs we will need the following lemmas.

\bs

\begin{lem1}
Let $\Th(z)$ be a meromorphic inner function with $1-\Th(t) \notin
L^2(\R)$ and let $\sigma$ be the corresponding Clark measure. If
$N[\Tb S^{2a}] \neq 0$ for some $a>0$, then
for any $\e>0$
there exists $h \in
L^{2}(\sigma)$ such that $$ \lim_{y \to \pm\infty }e^{xy} \int{
\frac{h(t)}{t-iy}d\sigma (t)} =0$$ for every $x \in (-a+\e,a-\e)$
and the measure $hd\sigma$ has finite total variation.
\end{lem1}

\begin{proof} The idea of the proof is truly simple: If the Toeplitz
kernel from the statement is non-trivial then $K_\Theta$
contains a function divisible by $S^{2a}$. The desired measure $hd\sigma$
is then obtained from the Clark representation of that function.
The details are as follows.

Let
$$b(z):=\frac{z-i}{z+i}.$$
Since $N[\Tb S^{2a}] \neq 0$, $N[\Tb S^{2a-2\e}b] \neq 0$ for any
$\e>0$. Hence there exists a non-zero
$$f \in  N[\Tb S^{2a-2\e}b]\subset
H^{2}(\C _{+}).$$
 Then $S^{2a-2\e}bf \in
K_\Th.$ Define $h:=S^{a-\e} bf/(z-i)$. Clearly
 $h$ belongs to $K_\Th$, and therefore $$ h(z)=\frac{1-\Th (z)}{2\pi i} \int{
\frac{h(t)}{t-z} d\sigma (t)} $$
where $\sigma$ is the Clark measure of $\Th$. In particular, for $x<a-\e,$  $$
\lim_{y \to \infty} e^{xy} \int{\frac{h(t)}{t-iy}d\sigma
(t)}=0$$  because $f(iy)\to 0$, since $f\in H^2(\C_{+})$,
and because the outer function $1-\Th(iy)$ cannot go to zero
exponentially fast.

Denote $g=\bar\Th h\in \bar H^2=H^2(\C_-)$. Then $h=\Th g$ in the lower half-plane. Note that
$g=S^{-a+\e} k$ where
$$k=\bar\Th S^{2a-2\e}bf/(z-i)\in \bar H^2=H^2(\C_-).$$

\no Hence
 for $x>-a+\e,$ $$ \lim_{y \to -\infty} e^{xy}
\int{\frac{h(t)}{t-iy}d\sigma (t)}=\lim_{y \to -\infty}2\pi i
\frac{k(iy)\Th (iy)}{e^{(a-\e+x)y}(1-\Th(iy))}=0.$$ The last equality
follows from the facts that $k(z) \in H^{2}(\C _{-})$ and
that $$\frac{1-\Th(iy)}{\Th(iy)}=\overline\Th(-iy)-1.$$

It is left to notice that $h(z)=l(z)/(z-i)$ where both
$l(t)=S^{2a-2\e}bf$ and $(z-i)^{-1}$ belong to $L^2(\sigma)$. Thus $h\in L^1(\sigma)$ and
 $hd\sigma$ has  finite total variation.
\end{proof}

\bs

The following Lemma is a well known fact whose proof we include
here for completeness.

\bs

\begin{lem2} Let $\mu$ be a measure with finite total variation.
Then the Fourier transform of $\mu$ vanishes on $[-a,a]$ if and
only if $$ \lim_{y \to \pm\infty }e^{xy}
\int{\frac{d\mu(t)}{t-iy}} =0,$$ for every $x \in [-a,a].$
\end{lem2}

\begin{proof} Suppose that $\int e^{ixt}d\mu(t)=0$ for all $x \in
[-a,a]$. Then
$$e^{-ixz}\int_{-\infty}^{+\infty}{\frac{e^{ixt}-e^{ixz}}{i(t-z)}}=\int_{-\infty}^{+\infty}\int_{0}^{x}e^{iu(t-z)}du d\mu(t)=$$
$$=\int_{0}^{x}\int_{-\infty}^{+\infty}e^{iut}d\mu(t)e^{-iuz}du=0,$$
for every $x \in [-a,a]$ and $z \in \C$. Therefore,
$$\int{\frac{e^{ixt}-e^{ixz}}{t-z}d\mu(t)}=0$$ for every $x \in [-a,a]$. Obviously,
\begin{equation}\lim_{y \to \pm\infty}\int{\frac{e^{ixt}}{t-iy}d\mu(t)}=0\label{eqn1}\end{equation}
 and therefore $$ \lim_{y \to \pm\infty }e^{xy}
\int{\frac{d\mu(t)}{t-iy}} =0$$ for every $x \in [-a,a].$

Conversely, for $x \in [-a,a]$, define
$$H(z):=\int{\frac{e^{ixt}-e^{ixz}}{t-z}d\mu(t)}.$$
Then $H(z)$ is an entire function of Cartwright class. To show
that $H(z)$ is identically zero it suffices to check that $\lim_{y
\to \pm\infty}H(iy)=0.$
Recall that for $x \in [-a,a]$,
$$ \lim_{y \to
\pm\infty}\int{\frac{e^{ix(iy)}}{t-iy}d\mu(t)}=0.$$ Together with \eqref{eqn1} this implies
$H\equiv 0$.

Thus
 $$\int
{e^{ixt}d\mu(t)}=\lim_{y \to \infty} -iy
\int{\frac{e^{ixt}d\mu(t)}{t-iy}}=\lim_{y\to
\infty}-iye^{-xy}\int{\frac{d\mu(t)}{t-iy}} =0$$ for all $-a\leq x
\leq a$.
\end{proof}

\bs\subsection{Main proofs}\no
Now we are ready to prove our main theorems. We will do it in the reverse order.

\bs

\begin{proofC} The inequality $\TT(X)\leq \GG(X)$ follows from Lemma~1 and Lemma~2.
To prove the opposite inequality, let $\GG(X)=a$. Then for any
$\e>0$ there exists a non-zero complex measure of total variation
no greater than $1$ supported on $X$ whose Fourier transform
vanishes on $[0,a-\e]$. Consider the set of all such measures.
Since this set is closed, convex and contains non-zero elements,
by the Krein-Milman theorem it has an extreme point, a non-zero
measure $\nu$.  Similarly to the proof of Theorem~66 in
\cite{dBr}, we can show that the extremality of $\nu$ implies that
it is supported on a discrete subset of $X$. Let $\Th(z)$ be the
meromorphic inner function whose Clark measure is $|\nu|$.
Then $\{\Th=1\}\subset X$. It is left to notice that the function
$$f(z)=\frac{1-\Th (z)}{2\pi i} \int{
\frac{d\nu(t)}{t-z} } $$ belongs to $K_\Th$ and is divisible by
$S^{a-\e}$ (as follows, for instance, from the proof of Lemma~2).
Hence $f/S^{a-\e}\in N[\Tb S^{a-\e}]\neq 0$.
\end{proofC}

\bs

\begin{proofB} (i) By Theorem C it is enough to prove that $\TT(\L)\geq 2\pi D_*(\L)$. Suppose that $D_*(\L)=a/2\pi$.
By the second definition \eqref{id2} of $D_*$ the function
$$\phi(x)=-2\pi n_\L(x)+(a-\e)x$$ is almost decreasing
for any $\e>0$. Consider a meromorphic inner function $\Th$ with
$\{\Th=1\}=\L$ and bounded derivative on $\R$. Then $\arg \bar\Th
S^{a-\e}$ differs from $ \phi$ by a bounded function. Hence $\arg
\bar\Th S^{a-2\e}$ is almost decreasing. By Theorem~II and
Theorem~III
$$N[\Tb S^{a-3\e}]\neq 0.$$

\bs

(ii) Again we will prove that $\TT(X)\leq 2\pi D_*(X)$. If
$\TT(X)=a$ then for any $\e>0$ there exists a meromorphic inner
$\Th(z)$ such that $\Gamma :=\{\Th=1\}\subset X$ and
$$N[\Tb S^{a-\e}]\neq 0.$$
By Theorem II (and remark after it) this means that $\arg \Tb S^{a-2\e}$ is almost decreasing. Hence
$$-2\pi n_\Gamma(x) +(a-3\e)x$$
 is
almost decreasing. Since $\e$ is arbitrary, $D_*(X)\geq a/2\pi$.
\end{proofB}

\ms

\begin{remark} It was pointed out by the referee that a different proof of Theorem~B
can be obtained from  Theorems~66 and 67 in \cite{dBr} and a
theorem of Krein.

\end{remark}

\begin{proofA}
$(ii)\Leftrightarrow(iii)$ follows from Theorem B and $(ii)\Leftrightarrow(iv)$
from Theorem C.

$(i)\Rightarrow(iii)$ Assume $(iii)$ is not true, i.e. for every
meromorphic inner function $\Th(z)$ with $\{\Th=1\}=\Lambda$,
$N[\bar{\Th}S^{2c}]=0$ for every $c>0.$ In this case we will
construct a non-constant zero type entire function which is
bounded on $\Lambda$, which will mean that $\Lambda$ is not a
P\'olya set. Define a measure $\mu$ to be  the counting measure of
$\Lambda$.  Then clearly $\int{d\mu(t)/(1+t^2)}<\infty.$ By
Theorem~I there exists a short de~Branges space $B_E$ contained
isometrically in $L^2(\mu)$. First, let us show that $B_E$ cannot
contain a function of positive exponential type. If $E(z)$
has type zero then all functions in $B_E$ have type zero by a theorem of Kein. If the type of $E$ is positive, then by Theorem~I we can assume that $B_E$ is contained properly in $L^2(\mu)$.

Suppose that $F(z)\in B_E$ has positive type.
We can assume that $F(iy)$ grows exponentially
in $y$ as $y\to\infty$.
Since $B_E\neq
L^2(\mu)$, there exists $g\in L^2(\mu)$ with $\bar{g}\perp B_E$.
Then

$$0=\int \frac{F(t)-F(w)}{t-w}g(t)d\mu(t)=\int \frac{F(t)}{t-w}g(t)d\mu(t)-F(w)\int \frac{1}{t-w}g(t)d\mu(t)$$

\no for any $w\in\C$ and therefore

$$F(w)=\frac{\int \frac{F(t)}{t-w}g(t)d\mu(t)}{\int \frac{1}{t-w}g(t)d\mu(t)}.$$
Since $F(w)$ grows exponentially along $i\R_+$, the integral in the denominator must decay exponentially in $w$ along
$i\R_+$.
Thus the function

$$G(z):=\frac{1-\Th(z)}{2\pi i}\int \frac{1}{t-z}g(t)d\mu(t)$$
can be represented as $G(z)=S^c(z)h(z)$ for some nonzero $h(z)\in
H^2(\C_+)$ and $c>0$, and belongs to $K_\Theta$, where $\Theta(z)$ is the
inner function corresponding to the measure $\mu$. Hence $h\in
N[\bar\Theta S^c]$ and we have a contradiction.

Therefore any $F(z)\in B_E$ has zero type. It is left to notice
that

 $$|F(\lan)| \leq
\sqrt{\sum_m{|F(\lambda_m)|^2}}=\|F\|_{L^2(\mu)}<\infty, $$ which
means that $F(z)$ is bounded on $\Lambda.$

$(ii)\Rightarrow(i)$ This is Theorem~XI in~\cite{dBr1}. For
reader's convenience, we include de~Branges' proof. Let $F(z)$ be
a zero type entire function bounded on $\Lambda$ by some constant
$M>0$. For any integer $n \in \N$, $F^n(z)$ is also a zero type
function. Let $\mu$ be a nonzero measure with finite total
variation whose Fourier transform vanishes identically on $[-a,a]$
for some $a>0$. Then, by the proof of Lemma~2,
$\int{(e^{ixt}-e^{ixz})/(t-z)d\mu(t)}=0$ for every $x \in (-a,a).$
Define
$$H(z):=\int{\frac{F^n(t)-F^n(z)}{t-z}d\mu(t)}$$ for all $z \in \C$.
It is clear that $H(z)$ is an entire function of zero type. To
show that $H(z) \equiv 0$ it is enough to check that $H(iy) \to 0$
as $y \to \pm\infty$. This follows from
$$ \lim_{y \to \pm\infty} H(iy)=\lim_{y\to \pm\infty}\left[ \int{\frac{F^n(t)}{t-iy}d\mu(t)}
-F^n(iy)e^{xy}\int{\frac{e^{ixt}}{t-iy}d\mu(t)}\right]=0. $$
Therefore, $$\int{\frac{F^n(t)-F^n(z)}{t-z}d\mu(t)} \equiv 0.$$
Now,
$$\left|F(z)\left(\int{\frac{d\mu(t)}{t-z}}\right)^{1/n}\right|\leq
M\left(\frac{\|\mu\|}{|\Im{z}|}\right)^{1/n},$$ for every non real
$z$. Since this is true for all $n \in \N$, we have that $|F(z)|
\leq M$ for all non real $z \in \C$ for which $\int{d\mu(t)/(t-z)}
\neq 0$. Since $\mu$ is a non-zero measure, by continuity, $F(z)$
is bounded in the whole plane.
\end{proofA}

\bs


\begin{thebibliography}{24}



\bibitem{Beu}{\sc Beurling, A.} {\it On quasianalyticity and
general distributions,} Mimeographed lecture notes, Summer
institute, Stanford University (1961)

\bibitem{BM1}  {\sc Beurling, A., Malliavin, P.} {\it
On Fourier transforms of measures with compact support,}
Acta Math.  107 (1962), 291--302

\bibitem{BM2}  {\sc Beurling, A., Malliavin, P.} {\it
On the closure of characters and the zeros of entire functions,}
Acta Math.  118 (1967), 79-93

\bibitem{Boa}  {\sc Boas, R. P.} {\it Entire Functions,} Academic
Press, New York, 1954





\bibitem{dBr}  {\sc De Branges, L.} {\it Hilbert spaces of entire functions.} Prentice-Hall,
Englewood Cliffs, NJ, 1968

\bibitem{dBr1} {\sc De Branges, L.} {\it Some applications of spaces
of entire functions,} Canadian J. Math. 15 (1963), 563-583

\bibitem{dBr2} {\sc De Branges, L.} {\it Some Hilbert spaces of
entire functions II,} Trans. Amer. Math. Soc. 99 (1961), 118-152

\bibitem{dBr4} {\sc De Branges, L.} {\it Some Hilbert spaces of
entire functions IV,} Trans. Amer. Math. Soc. 105 (1962), 43-83

\bibitem{Ga}{\sc Garnett, J.  } {\it Bounded analytic functions.} Academic Press, New York, 1981


\bibitem{HJ} {\sc Havin, V.,  J\"oricke, B.} {\it The uncertainty principle in harmonic analysis.}
 Springer-Verlag, Berlin, 1994.

\bibitem{MH} {\sc   Havin, V., Mashreghi, J.}  {\it Admissible majorants for model subspaces of $H^2$; I. Slow winding of the generating inner function, II. Fast winding of the generating inner function,} Canad. J. Math. 55 (2003), 1231--1263, 1264--1301.







\bibitem{Koo1} {\sc  Koosis, P.} {\it The logarithmic integral, Vol. I.} Cambridge Univ. Press, Cambridge, 1988





\bibitem{Lv1} {\sc Levin, B.} {\it Lectures on entire functions}
AMS, Providence, RI, 1996


\bibitem{Lev}{\sc Levinson, N.} {\it Gap and density theorems,} AMS
Colloquium Publications, 26 (1940)


\bibitem{MP} {\sc  Makarov, N.,  Poltoratski, A.} {\it Meromorphic inner functions, Toeplitz kernels, and the uncertainty principle,} in {\it Perspectives in Analysis}, Springer Verlag, Berlin, 2005, 185--252

\bibitem{MP2} {\sc  Makarov, N.,  Poltoratski, A.} {\it  Beurling-Malliavin theory for Toeplitz kernels,}
preprint, arXiv:math/0702497




\bibitem{Ni2} {\sc Nikolskii, N. K.} {\it Treatise on the shift
operator},  Springer-Verlaag, Berlin (1986)

\bibitem{Ni}  {\sc  Nikolskii, N.} {\it Operators, functions, and systems: an easy reading, Vol. I \& II.} AMS, Providence, RI, 2002




\bibitem{P} {\sc Poltoratski, A.} {\it Spectral gaps for sets and measures,} Preprint

\bibitem{Pol} {\sc P\'olya, G.} {\it Jahresbericht der Deutchen
Mathematiker-Vereinigung}, Vol. 40 (1940), Problem 105





\bibitem{Val} {\sc Valiron, G.}  {\it Sur la formule d'interpolation de Lagrange,} Bull. Sci. Math.
49 (1925), 181-192, 203-224


\end{thebibliography}
\end{document}